\documentclass[journal,twoside,web]{ieeecolor}
\usepackage{lcsys}
\usepackage{cite}
\usepackage{amsmath,amssymb,amsfonts}
\usepackage{algorithmic}
\usepackage{graphicx}
\usepackage{textcomp}

\newcommand{\ignore}[1]{}

\newcommand{\norm}[1]{\left\Vert#1\right\Vert} 

\newcommand{\bbm}{\begin{bmatrix}}
\newcommand{\ebm}{\end{bmatrix}}
\newcommand{\bma}[1]{\left[\begin{array}{#1}}
\newcommand{\ema}{\end{array}\right]}

\DeclareMathAlphabet{\mbf}{OT1}{ptm}{b}{n}
\newcommand{\mbs}[1]{{\boldsymbol{#1}}}
\newcommand{\mbc}[1]{ \boldsymbol{\mathcal{#1}} } 

\newcommand{\mbsdot}[1]{{\dot{\boldsymbol{#1}}}}

\newcommand{\mbfdot}[1]{{\dot{\mbf{#1}}}}
\newcommand{\mbfbar}[1]{{\bar{\mbf{#1}}}}
\newcommand{\mbfhat}[1]{{\hat{\mbf{#1}}}}

\newcommand{\mbftilde}[1]{{\tilde{\mbf{#1}}}}

\def\fdotb{{\raisebox{-0.6ex}{ \kern0.2ex\raisebox{0.8ex}{\tiny $\hspace*{-1ex}\circ$}}}}

\def\fddotb{{\raisebox{-0.6ex}{ \kern0.2ex\raisebox{0.8ex}{\tiny $\hspace*{-1ex}\circ\circ$}}}}

\newcommand{\f}{\frac}

\newcommand{\trans}{{\ensuremath{\mathsf{T}}}} 
\newcommand{\beq}{\begin{equation}}
\newcommand{\eeq}{\end{equation}}
\newcommand{\bdis}{\begin{displaymath}}
\newcommand{\edis}{\end{displaymath}}
\newcommand{\beqarray}{\begin{eqnarray}}
\newcommand{\eeqarray}{\end{eqnarray}}
\newcommand{\beqarraynn}{\begin{eqnarray*}}
\newcommand{\eeqarraynn}{\end{eqnarray*}}
\newcommand{\eye}{\mbf{1}}

\def\BibTeX{{\rm B\kern-.05em{\sc i\kern-.025em b}\kern-.08em
    T\kern-.1667em\lower.7ex\hbox{E}\kern-.125emX}}
\begin{document}
%
%
%
%
%
%
%
\def \myJournal {IEEE Control Systems Letters}
\def \myDoi {10.1109/LCSYS.2022.3178785}
\def \myPaperSiteName {IEEE Xplore}
\def \myPaperSiteLink {https://ieeexplore.ieee.org/document/9784861}
\def \myYear {2022}
\def \myPaperCitation{F. Ahmed, L. Sobiesiak and J. R. Forbes, ``Model Predictive Control of a Tandem-Rotor Helicopter With a Nonuniformly Spaced Prediction Horizon,'' in \textit{IEEE Control Systems Letters}, vol. 6, pp. 2828-2833, 2022.}


\begin{figure*}[t]

\thispagestyle{empty}
\begin{center}
\begin{minipage}{6in}
\centering
This paper has been accepted for publication in \emph{\myJournal}. 
\vspace{1em}

This is the author's version of an article that has, or will be, published in this journal or conference. Changes were, or will be, made to this version by the publisher prior to publication.
\vspace{2em}

\begin{tabular}{rl}
DOI: & \myDoi\\
\myPaperSiteName: & \texttt{\myPaperSiteLink}
\end{tabular}

\vspace{2em}
Please cite this paper as:

\myPaperCitation

\vspace{15cm}
\copyright \myYear \hspace{4pt}IEEE. Personal use of this material is permitted. Permission from IEEE must be obtained for all other uses, in any current or future media, including reprinting/republishing this material for advertising or promotional purposes, creating new collective works, for resale or redistribution to servers or lists, or reuse of any copyrighted component of this work in other works.

\end{minipage}
\end{center}
\end{figure*}
\newpage
\clearpage
\pagenumbering{arabic} 
	
\title{Model Predictive Control of a Tandem-Rotor Helicopter with a Non-Uniformly Spaced Prediction Horizon}
\author{Faraaz Ahmed, Ludwik Sobiesiak, and James Richard Forbes
\thanks{This work is supported by NGC Aerospace as well as the Mitacs Accelerate and NSERC Discovery Grant programs.}
\thanks{Faraaz Ahmed (email: faraaz.ahmed@mail.mcgill.ca) and James Richard Forbes (email: james.richard.forbes@mcgill.ca) are with the Department of Engineering, McGill University, Montreal QC, Canada, H3A 0C3.}
\thanks{Ludwik Sobiesiak (email: ludwik.sobiesiak@ngcaerospace.com) is with NGC Aerospace Ltd., Sherbrooke QC, Canada, J1L 2T9.}}

\maketitle
\thispagestyle{empty}

\begin{abstract}
This paper considers model predictive control of a tandem-rotor helicopter. The error is formulated using the matrix Lie group $SE_2(3)$. A reference trajectory to a target is calculated using a quartic guidance law, leveraging the differentially flat properties of the system, and refined using a finite-horizon linear quadratic regulator. The nonlinear system is linearized about the reference trajectory enabling the formulation of a quadratic program with control input, attitude keep-in zone, and attitude error constraints. A non-uniformly spaced prediction horizon is leveraged to capture the multi-timescale dynamics while keeping the problem size tractable. Monte-Carlo simulations demonstrate robustness of the proposed control structure to initial conditions, model uncertainty, and environmental disturbances.
\end{abstract}                                                    

\vspace{-1em}
\section{Introduction}
\label{sec:introduction}
\IEEEPARstart{T}{he} number of applications for unmanned aerial vehicles (UAVs) are growing, now including delivery, search and rescue, surveillance, and inspection \cite{Shakhatreh2019, Cherif2021}. While quadrotor and traditional helicopter platforms are popular for these types of tasks, a tandem rotor helicopter offers several advantages, including a large center of mass range, and larger lift capacities with smaller rotors \cite{Johnson2013}.

Matrix Lie groups can be used to compactly and accurately represent vehicle attitude, pose, or extended pose, which in turn can be leveraged in state estimation and control problems \cite{Bonnabel2009, Barrau2017a, Diemer2015, Cohen2020c, Chang2020, Kalabic2017}. In particular, an invariant linear quadratic Gaussian controller defined on $SE(2)$ is used to control a simplified car in \cite{Diemer2015}, and invariant linear quadratic regulator (ILQR) using an error defined on $SE_2(3)$ is used to control a quadrotor in \cite{Cohen2020c}. The use of an invariant error definition, instead of a traditional multiplicative error definition \cite{Farrell2008}, along with a particular type of process model, results in Jacobians that are state-independent. This yields improved robustness to initial conditions for state-estimation and control. Model predictive control (MPC) on matrix Lie groups has been explored in \cite{Chang2020} and \cite{Kalabic2017} for spacecraft attitude control on $SO(3)$.

In this paper, an invariant error definition \cite{Bonnabel2009, Barrau2017a} is used to develop an MPC strategy for a tandem-rotor helicopter. The attitude, velocity, and position states are cast into an element of the matrix Lie group $SE_2(3)$ \cite{Barrau2017a}. The nonlinear dynamics are linearized about a reference trajectory, allowing the MPC optimization problem to be posed as a quadratic program (QP). The QP is subject to various linear constraints on the control inputs and states. The dynamics used for control design do not exactly fit the invariant framework. Nevertheless, the invariant approach to control is followed due to the straightforward linearization process and the reduced dependence of the Jacobians on attitude \cite{Diemer2015, Cohen2020c}.

Unlike \cite{Cohen2020c}, where the angular velocity is a control input to the plant, the proposed MPC algorithm is able to control and constrain the body force and torque directly. This is accomplished by augmenting the $SE_2(3)$ state with the vehicle angular momentum, eliminating the need for a separate inner-loop controller to generate torque commands. A challenge with this approach is that the resulting process model is multi-timescale, necessitating a small controller timestep and long MPC prediction horizon \cite{Zlotnik2017}. This combination is computationally burdensome and is remedied by introducing a non-uniformly spaced prediction horizon, much like \cite{Tan2016} and \cite{Brudigam2021}.

The contributions of this paper are as follows. First, the synthesis of an MPC strategy that employs a non-uniformly spaced prediction horizon is presented. The system is linearized using an augmented $SE_2(3)$ error definition along a reference trajectory. The reference trajectory, which can be recomputed online, is generated using a combination of a quartic polynomial \cite{Lafontaine2004} and the solution to a finite-horizon LQR problem on $SE_2(3)$ \cite{Cohen2020c}. The MPC formulation features control input and state constraints. The second contribution is the inclusion of an attitude keep-in zone \cite{Walsh2017} and an $\ell_1$-norm \cite{Boyd2004} constraint on the attitude error, which together enforce attitude constraint satisfaction.
\vspace{-0.5em}
\section{Preliminaries} 
\label{sec:preliminaries}
\subsection{Kinematics and Reference Frames}
An inertial frame $\mathcal{F}_a$ is defined using the North-East-Down basis vectors \cite{Hughes2004}. Let point $w$ represent a point in $\mathcal{F}_a$ \cite{Bernstein2008}.  The frame $\mathcal{F}_b$ is defined by the set of orthonormal basis vectors $\underrightarrow{b}_1$, $\underrightarrow{b}_2$, and $\underrightarrow{b}_3$ that point out of the nose, right side, and underside of the fuselage respectively. Frame $\mathcal{F}_b$ is fixed to and rotates with the body of the helicopter. A physical vector $\underrightarrow{v}$ can be resolved in either $\mathcal{F}_a$ as $\mbf{v}_a$, or in $\mathcal{F}_b$ as $\mbf{v}_b$. The two are related by $\mbf{v}_a = \mbf{C}_{ab}\mbf{v}_b$, where $\mbf{C}_{ab} \in SO(3)$ is the direction cosine matrix (DCM) relating $\mathcal{F}_a$ and $\mathcal{F}_b$.
\vspace{-1em}
\subsection{Tandem Rotor Equations of Motion}
The tandem rotor helicopter is modeled as a rigid body subject to thrust, gravitational, and drag forces.  Let $z$ be a point collocated with the center of the mass of the helicopter. The kinematics are \cite{Hughes2004}
\begin{align}
	\mbfdot{C}_{ab} = \mbf{C}_{ab}\mbs{\omega}^{ba^\times}_b, \qquad
	\mbfdot{r}^{zw}_a = \mbf{v}^{zw/a}_a, \label{eq201}
\end{align}
where $\mbs{\omega}^{ba}_b$ is the angular velocity of $\mathcal{F}_b$ relative to $\mathcal{F}_a$ resolved in $\mathcal{F}_b$, $\mbf{r}^{zw}_a$ is the position of point $z$ relative to point $w$ resolved in $\mathcal{F}_a$, and $\mbf{v}^{zw/a}_a$ is the velocity of point $z$ relative to point $w$ with respect to $\mathcal{F}_a$, resolved in $\mathcal{F}_a$.  The cross operator $(\cdot)^\times: \mathbb{R}^3 \rightarrow \mathfrak{so}(3)$ is defined such that $\mbf{u}^\times \mbf{v} = -\mbf{v}^\times \mbf{u}$. 
The dynamics are \cite{Hughes2004}
\begin{align}
	m_\mathcal{B}\mbfdot{v}^{zw/a}_a = \mbf{C}_{ab}\mbf{f}^{\mathcal{B}z}_b, \quad
	\mbf{J}^{\mathcal{B}z}_b \mbsdot{\omega}^{ba}_b = \mbf{m}^{\mathcal{B}z}_b - \mbs{\omega}^{ba^\times}_b \mbf{J}^{\mathcal{B}z}_b \mbs{\omega}^{ba}_b,\label{eq202}
\end{align}
where $m_\mathcal{B}$ is the helicopter's mass, and $\mbf{J}^{\mathcal{B}z}_b$ is the helicopter's second moment of mass resolved in $\mathcal{F}_b$. The forces acting on the helicopter are $\mbf{f}^{\mathcal{B}z}_b = \mbf{f}^p_b + \mbf{C}_{ab}^\trans\mbf{f}^a_a + \mbf{C}_{ab}^\trans\mbf{f}^g_a$, where $\mbf{f}^p_b = [0 \; 0 \; f]^\trans$ is the propulsion force, $f$ is the total thrust force from the rotors, $\mbf{f}^a_a = \mbf{C}_{ab}\mbf{D}\mbf{C}_{ab}^\trans\mbf{v}^{zw/a}_a$ is the aerodynamic drag force, $\mbf{D} = \text{diag}(d_x, d_y, d_z)$ is a constant matrix composed of rotor drag coefficients \cite{Faessler2018}, and $\mbf{f}^g_a = [0 \; 0 \; m_\mathcal{B}g]^\trans$ is the gravitational force, where $g = 9.81 \text{ m/s}^2$. The torques acting on the helicopter are $\mbf{m}^{\mathcal{B}z}_b = \mbf{m}_b + \mbf{m}^a_b$, where $\mbf{m}_b$ is the total control moment from the rotors, and $\mbf{m}^a_b = -\mbf{E}\mbf{C}_{ab}^\trans\mbf{v}^{zw/a}_a - \mbf{F}\mbs{\omega}^{ba}_b$ is the parasitic torque from the rotor drag, where $\mbf{E}$ and $\mbf{F}$ are constant drag matrices \cite{Faessler2018}.
\vspace{-1em}
\subsection{Matrix Lie Groups}
Let $G$ denote a matrix Lie group and let $\mathfrak{g}$ denote the matrix Lie algebra associated with $G$ \cite{Barrau2015}. An element of $\mathfrak{g}$ can be mapped to $G$ using the matrix exponential, $\exp(\cdot): \mathfrak{g} \rightarrow G$, and the inverse operation is achieved using the matrix natural logarithm, $\log(\cdot): G \rightarrow \mathfrak{g}$. The matrix Lie algebra is mapped to a $k$ dimensional column matrix using the linear operator $(\cdot)^\vee: \mathfrak{g} \rightarrow \mathbb{R}^k$, and the inverse operation is performed using the operator $(\cdot)^\wedge: \mathbb{R}^k \rightarrow \mathfrak{g}$. For small $\delta\mbs{\xi} \in \mathbb{R}^k$, $\exp(\delta\mbs{\xi}^\wedge) \approx \eye + \delta\mbs{\xi}^\wedge$.

The group of double direct isometries, denoted as $SE_2(3)$ \cite{Barrau2015},  can be used to represent the attitude $\mbf{C}_{ab}$, velocity $\mbf{v}^{zw/a}_a$, and position $\mbf{r}^{zw}_a$ as
\begin{align}
	\mbf{X} = \bma{ccc} \mbf{C}_{ab} & \mbf{v}^{zw/a}_a & \mbf{r}^{zw}_a \\ \mbf{0} & 1 & 0 \\ \mbf{0} & 0 & 1 \ema \in SE_2(3). \label{eq203}
\end{align}
Details about the corresponding Lie algebra, $\mathfrak{se}_2(3)$, and the exponential map from $\mathfrak{se}_2(3) \rightarrow SE_2(3)$ are found in \cite{Barrau2015}.

\section{Control} 
\label{sec:control}
\subsection{Control Objective}
The objective of the controller is to generate actuator commands that allow the vehicle to follow a reference trajectory. Denote the desired reference frame $\mathcal{F}_r$. The reference attitude, velocity, and position trajectories are $\mbf{C}_{ar}$, $\mbf{v}^{z_rw/a}_a$, and $\mbf{r}^{z_rw}_a$, respectively. The reference states at timestep $k$ are written in terms of an $SE_2(3)$ element, $\mbf{X}^r_k$, using \eqref{eq203}. The tracking error is defined using a left-invariant error \cite{Barrau2017a}
\begin{align}
	\delta\mbf{X}_k = \mbf{X}^{r^{-1}}_k\mbf{X}_k = \bma{ccc} \delta\mbf{C}_k & \delta\mbf{v}_k & \delta\mbf{r}_k \\ \mbf{0} & 1 & 0 \\ \mbf{0} & 0 & 1 \ema \in SE_2(3), \label{eq411}
\end{align}
where $\delta\mbf{C}_k = \mbf{C}_{ar_k}^\trans\mbf{C}_{ab_k}$, $\delta\mbf{v}_k = \mbf{C}_{ar_k}^\trans(\mbf{v}^{zw/a}_{a_k} - \mbf{v}^{z_rw/a}_{a_k})$, and $\delta\mbf{r}_k = \mbf{C}_{ar_k}^\trans(\mbf{r}^{zw}_{a_k} - \mbf{r}^{z_rw}_{a_k})$. The tracking error is expressed in terms of the Lie algebra as $\delta \mbf{X}_k = \exp ( \delta \mbs{\xi}_k^\wedge)$ where $\delta \mbs{\xi}_k = [\delta \mbs{\xi}_k^{\phi ^\trans} \; \delta \mbs{\xi}_k^{v ^\trans} \; \delta \mbs{\xi}_k^{r ^\trans}]^\trans$. To incorporate the rotational dynamics, the state is augmented with the angular momentum tracking error, $\delta\mbf{h}_k$, defined as
\begin{align}
	\delta\mbf{h}_k = \delta\mbf{C}_k\mbf{h}^{\mathcal{B}z/a}_{b_k} - \mbf{h}^{\mathcal{B}z_r/a}_{r_k}, \label{eq416}
\end{align}
where $\mbf{h}^{\mathcal{B}z_r/a}_{r_k}$ is the reference angular momentum and $\mbf{h}^{\mathcal{B}z/a}_{b_k}$ is the true angular momentum. Therefore, the full state is $\delta\mbf{x}_k = [\delta\mbs{\xi}^{\phi^\trans}_k \; \delta\mbs{\xi}^{v^\trans}_k \; \delta\mbs{\xi}^{r^\trans}_k \; \delta\mbf{h}^\trans_k]^\trans$. The control objective is to drive the tracking error to zero such that $\delta\mbf{x}_k = \mbf{0}$.
\vspace{-1em}
\subsection{Control Inputs}
Unlike \cite{Cohen2020c}, where angular velocity tracking is achieved using a lower level inner-loop controller, the proposed controller outputs thrust force and torque commands directly
\begin{align}
	\mbf{u}_k = \bma{cc} f_k & \mbf{m}^{\trans}_{b_k} \ema^\trans \in \mathbb{R}^4. \label{eq418}
\end{align} 
This allows input constraints to be enforced at the torque command level.
\vspace{-1em}
\subsection{Overview of the Guidance and Control Structure }
\begin{figure}[t]
	\centering
	\centerline{\includegraphics[width=\columnwidth]{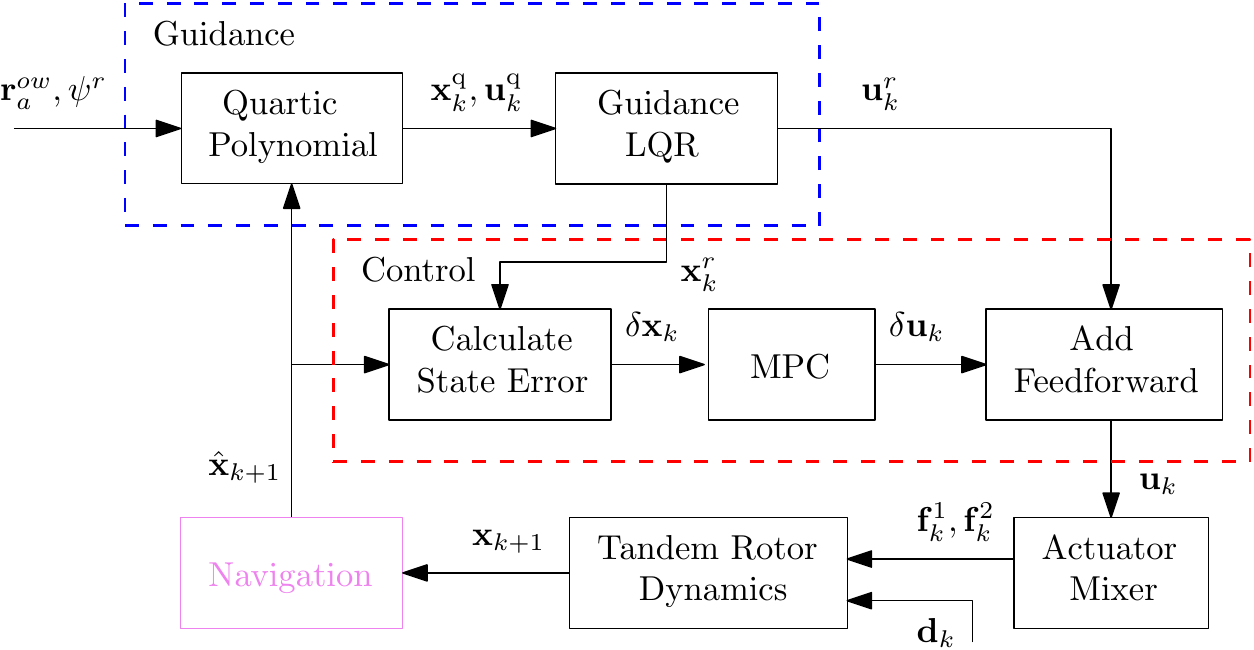}}
	\caption{Proposed guidance and control structure.}
	\label{fig1a}
	\vspace{-1em}
\end{figure}

Consider the guidance and control structure shown in Fig.~\ref{fig1a}. By using the control inputs in \eqref{eq418}, the dynamics of the tandem rotor helicopter are differentially flat \cite{Faessler2018}. A set of flat outputs, equal to the number of inputs, exists such that all of the system states and control inputs can be represented by the flat outputs and their derivatives \cite{Sferrazza2016}. For the tandem rotor helicopter, the flat outputs are $\mbs{\sigma}_k = [ \mbf{r}^{z_rw^\trans}_{a_k} \; \psi^r_k]^\trans$, where $\psi^r_k$ is the desired heading. 

A quartic guidance law is used to generate a coarse reference trajectory from the current position to the target position, $\mbf{r}^{z_fw}_a$, in terms of the flat outputs. The differentially flat property of the dynamics is used to generate the reference states, $\mbf{x}^q_k$, and control inputs, $\mbf{u}^q_k$. The reference trajectory is refined using finite-horizon LQR to create smooth trajectories for the states, $\mbf{x}^r_k$, and control inputs, $\mbf{u}^r_k$. 

The state error, $\delta\mbf{x}_k$, is calculated using $\eqref{eq411}$ and $\eqref{eq416}$. The MPC algorithm operates on the state error to produce the feedback control input, $\delta\mbf{u}_k$, which is added to the feedforward reference control input, $\mbf{u}^r_k$, to produce the total control input
\begin{align}
	\mbf{u}_k = \bma{c} f_k \\ \mbf{m}_{b_k} \ema = \bma{c} f^{r}_k + \delta f_k \\
	\delta\mbf{C}^\trans_k \mbf{m}^{r}_{r_k} + \delta\mbf{m}_{b_k} \ema. \label{eq420}
\end{align}
The reference torque command, $\mbf{m}^{r}_{r_k}$, is resolved in $\mathcal{F}_r$. Before being combined with $\delta\mbf{m}_{b_k}$, it must be resolved in $\mathcal{F}_b$ through multiplication by $\delta\mbf{C}^\trans_k$.

An actuator mixer is used to map the control inputs from \eqref{eq420} to front and rear rotor force components, $\mbf{f}^1_k$ and $\mbf{f}^2_k$, respectively. The tandem rotor helicopter is over-actuated, therefore it is assumed that the $\underrightarrow{b}_2$, and $\underrightarrow{b}_3$ rotor force components are used for control, and the $\underrightarrow{b}_1$ component is strictly used for trimming. The mapping from the control inputs, $\delta\mbf{u}_k$, to the rotor force inputs $\mbf{f}^{1,2}_k$ is as shown in \cite{Lee2010}.

In practice a navigation loop will provide state estimates, $\mbfhat{x}_k$, to the guidance and control algorithms. Herein it is assumed the state estimates are accurate such that $\mbf{x}_k = \mbfhat{x}_k$.

The guidance generates an initial reference trajectory at $t = 0$. If the current trajectory becomes invalid, possibly due to a large disturbance, a new trajectory is planned.

\vspace{-1em}
\subsection{Linearization of Dynamics}
Consider the following first-order approximations, valid for small $\delta\mbs{\xi}$ \cite{Hartley2019}, 
\begin{subequations}\label{eq439}
	\begin{align}
		\delta\mbf{C} &= \exp(\delta\mbs{\xi}^{\phi^\times}) \approx \eye + \delta\mbs{\xi}^{\phi^\times},\\
		\delta\mbf{v} &= \mbftilde{J}(\delta\mbs{\xi}^\phi)\delta\mbs{\xi}^v \approx \delta\mbs{\xi}^v,\\
		\delta\mbf{r} &= \mbftilde{J}(\delta\mbs{\xi}^\phi)\delta\mbs{\xi}^r \approx \delta\mbs{\xi}^r,
	\end{align}
\end{subequations}
where $\mbftilde{J}(\cdot)$ is the $SO(3)$ left Jacobian \cite{Hartley2019}. Using \eqref{eq439}, and the error definitions from \eqref{eq411}, \eqref{eq416}, and \eqref{eq420}, the continuous-time equations of motion are linearized about the reference trajectory yielding $\delta\mbfdot{x} = \mbf{A}\delta\mbf{x} + \mbf{B}\delta\mbf{u}$, where
\begin{align}
	\mbf{A} = \bma{ccccc} \mbf{0} & \mbf{0} & \mbf{0} & \mbf{A}_{14} \\ \mbf{A}_{21} & \mbf{A}_{22} & \mbf{0} & \mbf{0} \\ \mbf{0} & \eye & \mbf{A}_{33} & \mbf{0} \\ \mbf{0} & \mbf{0} & \mbf{0} & \mbf{A}_{44} \ema,\;
	\mbf{B} = \bma{cc} \mbf{0} & \mbf{0} \\ \mbf{B}_{21} & \mbf{0} \\ \mbf{0} & \mbf{0} \\ \mbf{0} & \eye \ema, \label{eq440}
\end{align}
$\mbf{A}_{14} = \mbf{J}^{\mathcal{B}z^{-1}}_b$, $\mbf{A}_{22} = -\mbs{\omega}^{ra^\times}_r - 1/m_\mathcal{B}\mbf{D}$, $\mbf{A}_{33} = -\mbs{\omega}^{ra^\times}_r$, $\mbf{A}_{44} = -\mbs{\omega}^{ra^\times}_r$, $\mbf{B}_{21} = -1/m_\mathcal{B}\eye_3$, and 
\begin{equation}
	\begin{aligned}
		\mbf{A}_{21} &= \f{1}{m_\mathcal{B}}\left(\left(\mbf{D}\mbf{C}_{ar}^\trans\mbf{v}^{z_rw/a}_a\right)^\times \right.\\ &\qquad\qquad \left. - \mbf{D}\left(\mbf{C}_{ar}^\trans\mbf{v}^{z_rw/a}_a\right)^\times + \left(f^r \eye_3\right)^\times \right).
	\end{aligned}
\end{equation}
The Jacobians from \eqref{eq440} only depend on reference quantities. Moreover, notice that only $\mbf{A}_{21}$ depends on $\mbf{C}_{ar}$, and when $\mbf{D} = \mbf{0}$, the Jacobians depend only on $f^r$ and $\mbs{\omega}^{ra}_r$.

The continuous-time linearized system is then discretized \cite{Farrell2008} yielding
\begin{align}
	\delta\mbf{x}_{k+1} = \mbf{A}_k \delta\mbf{x}_k + \mbf{B}_k \delta\mbf{u}_k. \label{eq422}
\end{align}
\vspace{-1em}
\subsection{Finite Horizon MPC for Linear Time Varying Systems}
For the discrete-time-linearized system \eqref{eq422}, the state predictions are $\delta\mbf{x}_{i|k} = \mbc{A}_i(k)\delta\mbf{x}_k + \mbc{C}_i(k)\delta\mbs{\mu}_k,\quad i = 0,\ldots,N,$ where $\delta\mbs{\mu}_k~=~[\delta\mbf{u}_{0|k}^\trans \; \cdots \; \delta\mbf{u}_{N-1|k}^\trans]^\trans$ is the predicted input sequence, $\delta\mbf{x}_{i|k}$ is the state vector at time $k+i$ predicted at time $k$, $N$ is the prediction horizon length, and the time-varying state transition matrices are
\begin{equation}\label{eq423}
	\begin{aligned}
		\mbc{A}_i(k) &= \overset{\curvearrowleft}{\prod^{i-1}_{j=0}} \mbf{A}_{k+j},\\
		\mbc{C}_i(k) &= \Bigg[\begin{matrix}
			\left(\overset{\curvearrowleft}{\prod^{i-1}_{j=1}} \mbf{A}_{k+j}\right)\mbf{B}_k & \left(\overset{\curvearrowleft}{\prod^{i-1}_{j=2}} \mbf{A}_{k+j}\right)\mbf{B}_{k+1}\end{matrix}\Bigg.\\ &\qquad\qquad
		\Bigg.\begin{matrix} \cdots & \mbf{B}_{k+i-1} & \mbf{0} & \cdots & \mbf{0} \end{matrix}\Bigg].
	\end{aligned}
\end{equation}
Note $\overset{\curvearrowleft}{\prod}$ indicates that successive terms in the sequence are left-multiplied. The cost matrices can be written as
\begin{align}
	\mbc{S}_k = \bma{c} \mbc{C}_1(k) \\ \vdots \\ \mbc{C}_N(k) \ema,\quad
	\mbc{M}_k = \bma{c} \mbc{A}_1(k) \\ \vdots \\ \mbc{A}_N(k) \ema.
\end{align}
Therefore, the predicted state sequence over the prediction horizon, $\delta\mbs{\chi}_k = [\delta\mbf{x}_{1|k}^\trans \; \cdots \; \delta\mbf{x}_{N|k}^\trans]^\trans$, can be written as
\begin{align}
	\delta\mbs{\chi}_k = \mbc{S}_k\delta\mbs{\mu}_k + \mbc{M}_k\delta\mbf{x}_k. \label{eq424}
\end{align}
The cost function to be optimized is
\begin{multline}
	J(\delta\mbs{\mu}_k) = \delta\mbf{x}^\trans_{N|k} \mbf{P} \delta\mbf{x}_{N|k} \\ + \sum^{N-1}_{i=0} \left( \delta\mbf{x}^\trans_{i|k} \mbf{Q}_k \delta\mbf{x}_{i|k} + \delta\mbf{u}^\trans_{i|k} \mbf{R}_k \delta\mbf{u}_{i|k} \right), \label{eq425}
\end{multline}
where $\mbf{Q} = \mbf{Q}^\trans \geq 0$ is the state penalty matrix, $\mbf{R} = \mbf{R}^\trans > 0$ is the control input penalty matrix, and $\mbf{P} = \mbf{P}^\trans \geq 0$ is the terminal state penalty matrix. Writing \eqref{eq425} in matrix form, the constrained optimization problem is expressed as a QP
\begin{equation}\label{eq431}
	\begin{aligned}
		&\underset{\delta\mbs{\mu}_k}{\text{min}}\quad J(\delta\mbs{\mu}_k) = \frac{1}{2}\delta\mbs{\mu}^\trans_k\mbf{H}_k\delta\mbs{\mu}_k + \mbf{F}^\trans_k\delta\mbs{\mu}_k,\\
		&\text{s.t.}\quad \mbf{G}_k\delta\mbs{\mu}_k \leq \mbf{W}_k + \mbf{T}_k\delta\mbf{x}_k,
	\end{aligned}
\end{equation}
where $\mbf{H}_k = \mbc{S}_k^\trans \mbfbar{Q} \mbc{S}_k + \mbfbar{R}$, $\mbf{F}_k = \mbc{S}_k^\trans \mbfbar{Q} \mbc{M}_k$, $\mbfbar{Q}~=~ \text{diag}\left(\mbf{Q},\ldots,\mbf{Q},\mbf{P}\right)$, and $\mbfbar{R} = \text{diag}\left(\mbf{R},\ldots,\mbf{R}\right)$. Solving the QP gives the optimal control input sequence $\delta\mbs{\mu}^\star_k$.
\vspace{-2em}
\subsection{Non-Uniform Prediction Horizon Timestep}
The process model given by \eqref{eq201} and \eqref{eq202} contains ``fast" attitude dynamics and ``slow" translational dynamics. A large prediction horizon is necessary to accurately predict the dynamics of this multi-timescale system. A non-uniformly spaced prediction horizon is implemented to increase the time span of the prediction horizon without increasing the number of optimization variables \cite{Tan2016}. The short-term horizon sampling time is kept small to resolve the ``fast" attitude dynamics, whereas the long-term horizon sampling time is increased to ensure the ``slow" translational dynamics are predicted sufficiently far into the future.

The penalty matrices $\mbf{Q}$ and $\mbf{R}$ are modified to account for the different sampling time at each segment in the horizon such that $\mbf{Q}_i = \left(\Delta t_i/\Delta t_1\right)\mbf{Q}$, and $\mbf{R}_i = \left(\Delta t_i/\Delta t_1\right)\mbf{R}$ \cite{Brudigam2021}, where $\Delta t_i$ is the sampling time of the individual horizon segment.
\vspace{-1em}
\subsection{State and Input Constraints}
An advantage of MPC is the ability to explicitly embed state and control input constraints in the optimization problem. The desired constraints must be defined in terms of the optimization variables and written as linear inequalities to be included in the QP. 

The attitude of the helicopter is constrained using a keep-in zone \cite{Walsh2017}, defined as
\begin{align}
	\mbf{x}^{1^\trans}_b \mbf{C}_{ab}^\trans \mbf{y}^1_a \geq \cos(\alpha) - \epsilon_1, \label{eq442}
\end{align}
where $\epsilon_1$ is a slack variable. By setting $\mbf{x}^1_b = \mbf{y}^1_a = [0 \; 0 \; 1]^\trans$, the roll and pitch angles are constrained simultaneously by $\alpha$. To incorporate this constraint in the QP, \eqref{eq442} is linearized using \eqref{eq411} and \eqref{eq439},
\begin{align}
	\mbf{x}^{1^\trans}_b \left(\eye - \delta\mbs{\xi}^{\phi^\times}\right) \mbf{C}_{ar}^\trans \mbf{y}^1_a \geq \cos(\alpha) - \epsilon_1.\label{eq445}
\end{align}
From the stacked state transition matrix \eqref{eq424}, the $i^\text{th}$ predicted attitude in the horizon can be written as
\begin{align}
	\delta\mbs{\xi}^\phi_{k|i} = \mbf{P}_i \left(\mbc{S}_k \delta\mbs{\mu}_k + \mbc{M}_k \delta\mbf{x}_k \right), \label{eq446}
\end{align}
where $\mbf{P}_i$ is a projection matrix. Substituting \eqref{eq446} into \eqref{eq445}, the linearized keep-in zone constraint is written in terms of the optimization variables as
\begin{multline}
	\bma{cc} -\mbf{x}^{1^\trans}_b \left( \mbf{C}_{ar_{k+i}}^\trans \mbf{y}^1_a \right)^\times \mbf{P}_i\mbc{S}_k & -\eye \ema \bma{c} \delta\mbs{\mu}_k \\ \epsilon_1 \ema \leq \\ \bma{c} -\cos(\alpha) + \mbf{x}^{1^\trans}_b \mbf{C}_{ar_{k+i}}^\trans \mbf{y}^1_a \ema \\ + \bma{c} \mbf{x}^{1^\trans}_b \left( \mbf{C}_{ar_{k+i}}^\trans \mbf{y}^1_a \right)^\times \mbf{P}_i \mbc{M}_k \ema \delta\mbf{x}_k.\label{eq447}
\end{multline}

Because \eqref{eq447} is only valid for small $\delta\mbs{\xi}^\phi$, an additional constraint is imposed on the size of the attitude error to ensure accuracy of \eqref{eq445}. While the $\ell_2$-norm is a logical choice, the $\ell_1$-norm is used here because it provides a conservative size constraint since $\norm{\mbf{x}}_1 \leq \sqrt{n}\norm{\mbf{x}}_2 , \forall \mbf{x} \in \mathbb{R}^{n_x}$. The attitude error constraint is written as
\begin{align}
	\norm{\delta\mbs{\xi}^\phi}_1 \leq \gamma + \epsilon_2, \label{eq448}
\end{align}
where $\epsilon_2$ is a slack variable. Equation \eqref{eq448} is equivalent to the linear inequalities \cite{Boyd2004}
\begin{align}
	\sum^3_{i=1} z_i \leq \gamma + \epsilon_2,\quad -z_i \leq \delta\mbs{\xi}^\phi_i \leq z_i, \label{eq449}
\end{align}
where $z_i, \; i = 1,2,3$ are optimization variables. Substituting \eqref{eq446} into \eqref{eq449}, the $\ell_1$ constraint is written in terms of the optimization variables as
\begin{multline}
	\bma{ccc} \mbf{0} & \eye^{1\times 3} & -1 \\ \mbf{P}_i\mbc{S}_k & -\eye & 0 \\ -\mbf{P}_i\mbc{S}_k & -\eye & 0 \ema \bma{c} \delta\mbs{\mu}_k \\ \mbf{z} \\ \epsilon_2 \ema \leq \\ \bma{c} \gamma \\ \mbf{0} \\ \mbf{0} \ema + \bma{c} 0 \\ -\mbf{P}_i\mbc{M}_k \\ \mbf{P}_i\mbc{M}_k \ema \delta\mbf{x}_k. \label{eq450}
\end{multline}
To define control input constraints, first a maximum, $\mbf{u}_\text{max}~=~[f_\text{max} \; \mbf{m}^{\trans}_{b_\text{max}}]^\trans$, and minimum, $\mbf{u}_\text{min}~=~[f_\text{min} \; \mbf{m}^{\trans}_{b_\text{min}}]^\trans$, total control effort is defined. The MPC algorithm operates on $\delta\mbf{u}_k$, therefore the control input constraints are $\delta\mbf{u}_{\text{max}_k} = \mbf{u}_{\text{max}} - \mbf{u}^r_k$, and $\delta\mbf{u}_{\text{min}_k} = \mbf{u}_{\text{min}} - \mbf{u}^r_k$. In matrix form, the input constraints over the control horizon are
\begin{align}
	\bma{c} \eye \\ -\eye \ema \delta\mbs{\mu}_k \leq \bma{c} \delta\mbs{\mu}_{\text{max}_k} \\ 
	-\delta\mbs{\mu}_{\text{min}_k} \ema. \label{eq451}
\end{align}
The constraints \eqref{eq447}, \eqref{eq450}, and \eqref{eq451} can then be included in the QP given by \eqref{eq431}.
\vspace{-1em}
\subsection{Disturbance Estimation}
To improve guidance, an estimate of the disturbances acting on the vehicle, particularly a near-constant wind, is needed. Assuming the states are accurately estimated, then the disturbance at the previous timestep can be approximated as $\mbf{d}_{k-1} = \mbf{x}_k - \mbf{f}(\mbf{x}_{k-1},\mbf{u}_{k-1})$, where $\mbf{f}(\mbf{x}_k,\mbf{u}_k)$ represents modeled dynamics from \eqref{eq201} and \eqref{eq202}. A simple moving average filter is used to estimate the current disturbance. When the trajectory is replanned, the controller passes the disturbance estimate to the LQR guidance to improve the accuracy of the new reference trajectory. Without this simple disturbance estimator in the guidance, the tracking performance is poor.
\vspace{-1em}
\subsection{Reference Trajectory Generation}
A quartic guidance law is used to provide a minimum-time position and velocity trajectory from the current position to a target location, as shown in \cite{Lafontaine2004}. Because the flat outputs are the position and heading, the remaining state trajectories, $\mbf{x}^q_k$, and control input trajectories, $\mbf{u}^q_k$, can be found as a function of the quartic polynomials and a specified heading, $\psi^r_k$ \cite{Faessler2018}.

The quartic trajectory is used to warm start a standard discrete-time, finite-horizon LQR problem \cite{Stengel1994}. The system dynamics from \eqref{eq440} are linearized about $\mbf{x}^q_k$ and $\mbf{u}^q_k$, and discretized. The optimal gain sequence, $\mbf{K}_\text{LQR}$, is then found using the approach from \cite{Cohen2020c}. The sequence of reference control inputs, $\mbf{u}^r_k$, and states, $\mbf{x}^r_k$, are generated by propagating the closed-loop dynamics, $\mbf{x}^r_{k+1} = \mbf{f}(\mbf{x}^r_k,\mbf{u}^r_k) + \mbfhat{d}_k$, where $\mbfhat{d}_k$ is the disturbance estimate from the previous control step.

\section{Simulation Results} 
\label{sec:simulations}
A series of simulations are performed using a tandem-rotor helicopter model and the equations of motion from \eqref{eq201} and \eqref{eq202}. The mass and second moment of mass are, $m_\mathcal{B} = 218 \text{ kg}$ and $\mbf{J}^{\mathcal{B}z}_b = \text{diag}(26.8, 97.6, 87.2) \text{ kg}\cdot\text{m}^2$, respectively. For each simulation, the target position and velocity is $\mbf{r}^{z_fw}_a = \mbf{0}\text{ m}$, and $\mbf{v}^{z_fw/a}_a = \mbf{0}\text{ m/s}$, respectively,  and the reference heading is $\psi^r = 0 \text{ rad}$.  The nominal initial position and velocity is $\mbf{r}^{z_0w}_a = [-30 \; -5 \; -20]^\trans \text{ m}$, and $\mbf{v}^{z_0w/a}_a = [5 \;\; 0 \; -0.5]^\trans \text{ m/s}$, respectively. The initial attitude and angular velocity is $\mbf{C}_{ab_0} = \eye$, and $\mbs{\omega}^{b_0a}_b = \mbf{0} \text{ rad/s}$ respectively. The nominal wind speed is set to $\mbf{v}^{sw/a}_a = [0 \; -5 \;\; 0]^\trans \text{ m/s}$ and wind gusts, $\mbf{v}^{gw/a}_a$, are generated using the Dryden model \cite{Moorhouse1982}.

The control loop is run at $50$ Hz. The MPC prediction horizon is $N = 20$, while the control horizon is limited to $N_u = 10$ to reduce computational complexity \cite{Schwenzer2021}. The $\ell_1$ attitude error constraint is set to $\gamma = 0.1 \text{ rad}$, while the keep-in zone is set to $\alpha = 0.14 \text{ rad}$. The control input constraints are $0 \leq f \leq 3000$ N and $-200 \leq m_{b_i} \leq 200 \text{ N}\cdot\text{m}$. The reference trajectory is replanned if the $\ell_1$ constraint is active for longer than $0.4$~s.
\vspace{-1em}
\subsection{Demonstrating Constraints and Trajectory Replanning}
\begin{figure}[t]
	\centering
	\includegraphics[width=\columnwidth]{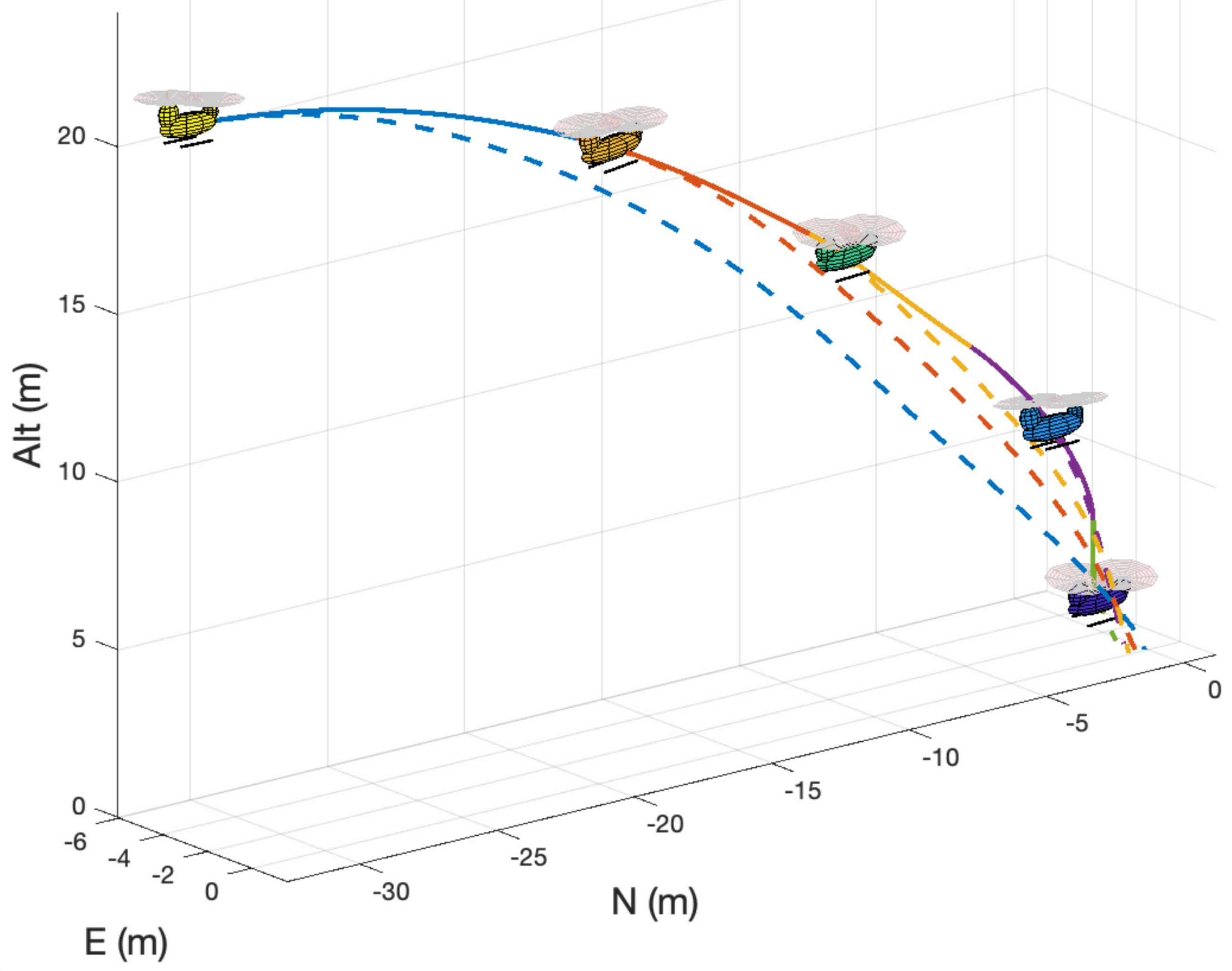}
	\vspace{-1em}
	\caption{Example of helicopter's progress along planned path. Dashed lines represent desired trajectory, solid lines represent actual trajectory. Changes in color represent replanned trajectories.}
	\label{fig1}
	\vspace{-0.5em}
\end{figure}

The results from a single simulation are highlighted to demonstrate advantages of the proposed control structure. The helicopter path from the starting position to the target is shown in Fig.~\ref{fig1}. Dashed lines show reference trajectories, while solid lines show actual trajectories. Changes in line color indicate a replanned trajectory, and occur when the attitude error exceeds the specified threshold. Initial tracking of the reference trajectory is good when the disturbance estimate from the controller is most accurate. As the disturbance evolves, the tracking performance suffers until the trajectory is replanned with an updated disturbance estimate.

\begin{figure}[t]
	\centering
	\includegraphics[width=\columnwidth]{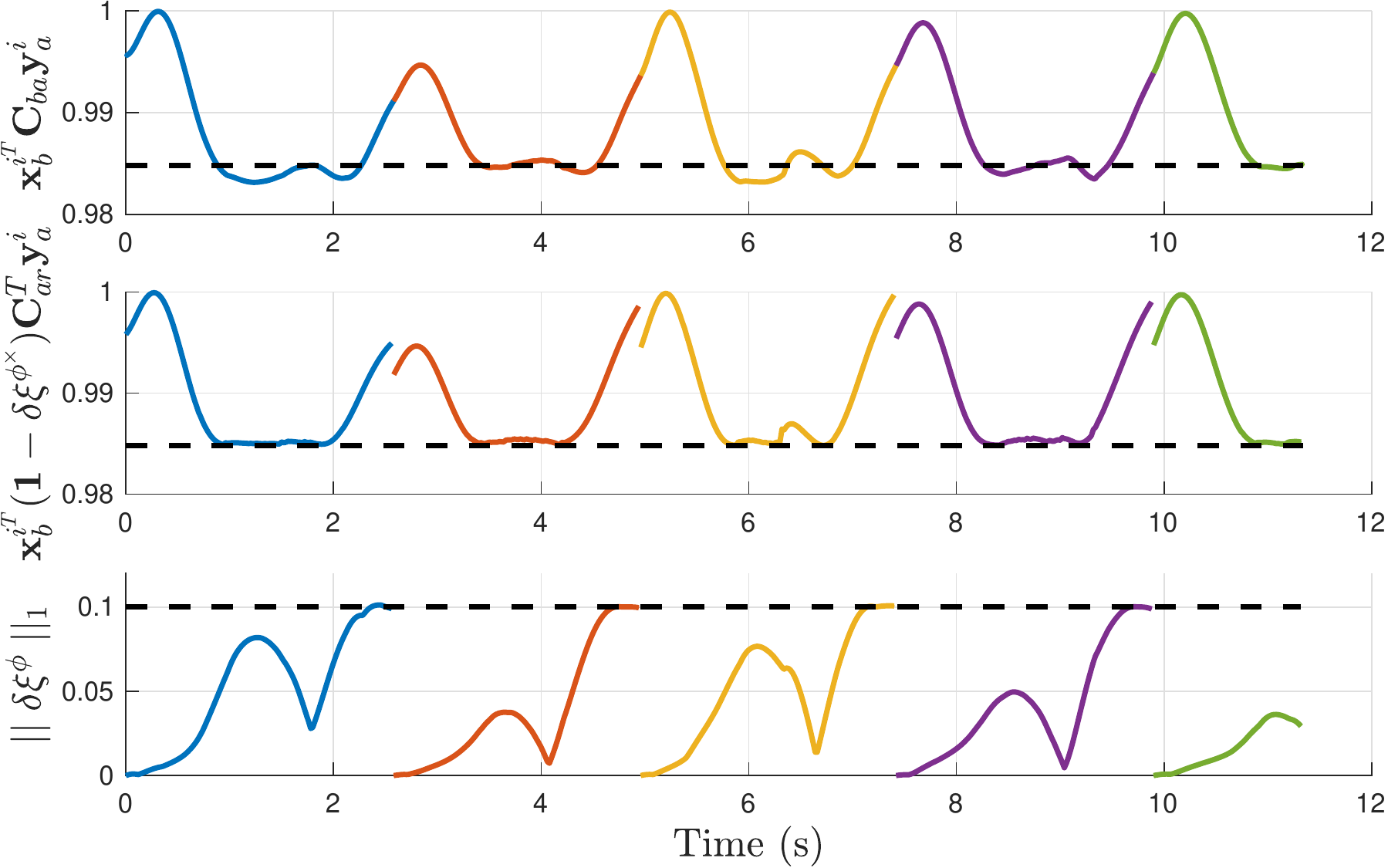}
	\vspace{-2em}
	\caption{Nonlinear (top) and linearized (middle) attitude keep-in zone constraints. The constraint is respected if the solid line is above the dashed line. The $\ell_1$ attitude error constraint (bottom) is respected if the solid line is below the dashed line.}
	\label{fig4}
	\vspace{-0.5em}
\end{figure}

The attitude keep-in zone and $\ell_1$ attitude error constraints are visualized in Fig.~\ref{fig4}. The linearized keep-in zone and attitude error constraints are respected throughout the entire simulation. During periods with larger attitude error, the linearized and nonlinear keep-in zone constraints diverge slightly. In some cases, the nonlinear constraint is marginally violated. This behavior is limited by the attitude error constraint, which maintains the validity of the keep-in zone by limiting the size of $\delta\mbs{\xi}^\phi$. The divergence of the nonlinear and linearized keep-in zones can be further limited by reducing $\gamma$. However, overly restricting the attitude error results in diminished tracking performance in the presence of large disturbances. It can be seen that the trajectory is replanned when the $\ell_1$ attitude error constraint becomes active. Once the trajectory is replanned, the attitude error immediately drops to zero since the new reference trajectory is planned from the current state. Although not shown here, the force and torque inputs are bound by the imposed constraints.
\vspace{-1em}
\subsection{Monte-Carlo Simulations}
Monte-Carlo simulations are performed to test the robustness of the proposed control structure to initial conditions, environmental disturbances, and model uncertainty. The initial state of the helicopter is randomized such that $\mbs{\phi}_0 = \mbf{0} + \mbf{w}_1 \text{ rad}$, $\mbf{r}^{z_0w}_a = [-30 \; -5 \; -20]^\trans + \mbf{w}_2 \text { m}$, $\mbf{v}^{z_0w/a}_a = [5 \; 0 \; 0.5]^\trans + \mbf{w}_3 \text{ m/s}$, and $\mbs{\omega}^{b_0a}_b = \mbf{0} + \mbf{w}_4 \text{ rad/s}$, where $\mbs{\phi}_0 = \log_{SO(3)}(\mbf{C}_{ab_0})$, $\mbf{w}_1 \sim \mathcal{N}(0,0.116^2\eye)$,  $\mbf{w}_2 \sim \mathcal{N}(0,\eye)$,  $\mbf{w}_3 \sim \mathcal{N}(0,0.333^2\eye)$,  and $\mbf{w}_4 \sim \mathcal{N}(0,0.029^2\eye)$. The nominal wind condition and Dryden gust model are randomized such that $\mbf{v}^{sw/a}_a = [0 \; -3 \; 0]^\trans + \mbf{w}_5$, and $W_0 \sim \mathcal{N}(10 \text{ m/s},1)$, where $\mbf{w}_5 \sim \mathcal{N}(0,1.667^2\eye)$, and $W_0$ is the low altitude intensity of the Dryden model. The estimated mass, $\hat{m}_\mathcal{B}$, and inertia matrix, $\mbfhat{J}^{\mathcal{B}z}_b$, are perturbed from their true values such that $\hat{m}_\mathcal{B} = m_\mathcal{B} + w_5$, and $\mbfhat{J}^{\mathcal{B}z}_{\hat{b}} = \mbf{C}^\trans_{b\hat{b}} \mbf{J}^{\mathcal{B}z}_b \mbf{C}_{b\hat{b}}$, 
where $w_5 \sim \mathcal{N}(0,10^2)$, and $\mbf{C}_{b\hat{b}} = \exp_{SO(3)}(\mbs{\phi})$ is a perturbation DCM, where $\mbs{\phi} \sim \mathcal{N}(0,0.044^2\eye)$.

\begin{figure}[!h]
	\centering
	\includegraphics[width=\columnwidth]{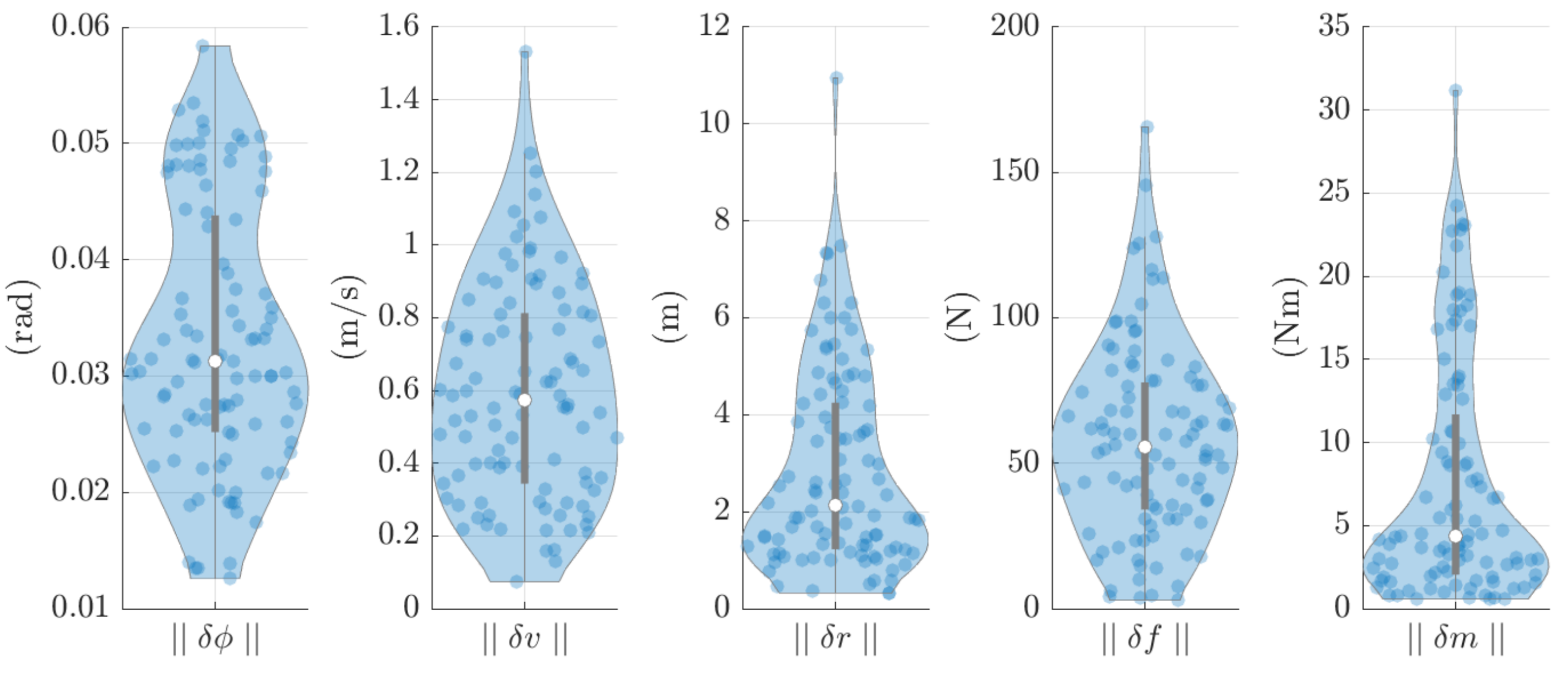}
	\vspace{-2em}
	\caption{State and input RMSE values across $100$ Monte-Carlo simulations with random initial conditions, wind gusts, and model uncertainty.}
	\label{fig5}
	\vspace{-1.5em}
\end{figure}

\begin{figure}[h]
	\vspace{-0.5em}
	\centering
	\includegraphics[width=\columnwidth]{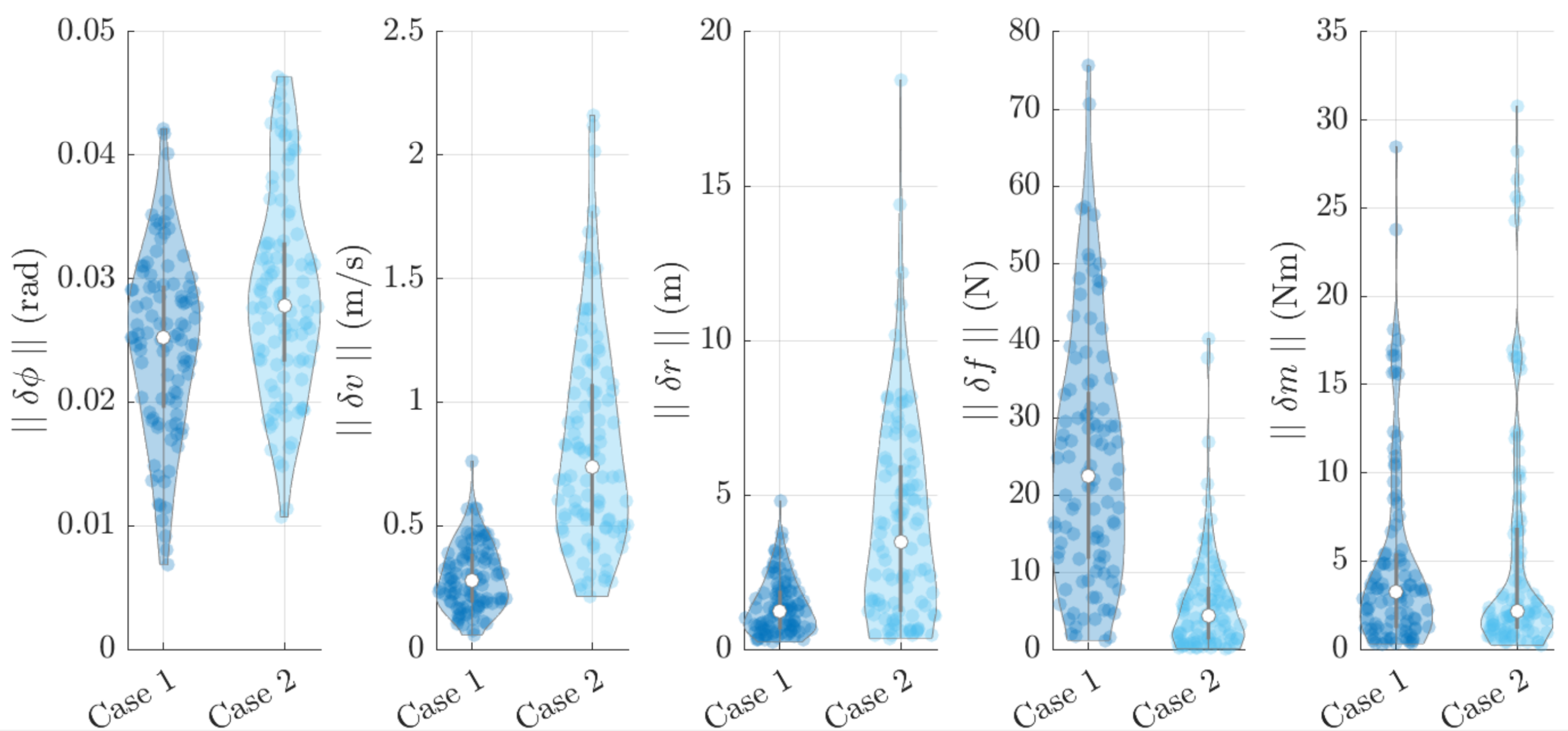}
	\vspace{-1.5em}
	\caption{State and input RMSE values across $100$ Monte-Carlo simulations comparing non-uniform (Case 1) and uniform (Case 2) prediction horizons.}
	\label{fig7}
	\vspace{-1em}
\end{figure}

The distribution of root-mean-square error (RMSE) results for $100$ Monte-Carlo runs is shown in Fig.~\ref{fig5} in the form of violin plots. In each case, the helicopter is able to reach the target position. The large distribution on the control input, particularly thrust, is due to the applied mass uncertainty. With reduced mass uncertainty, the variance in $\delta f$ and $\delta \mbf{m}_b$ is smaller, as seen in Fig.~\ref{fig7}. 

An additional set of Monte-Carlo simulations are run to demonstrate the advantage of the non-uniform prediction horizon. In Case $1$, the proposed non-uniform timestep is used, while Case $2$ features a uniform fixed-timestep equal to the controller timestep. The prediction horizon parameters for both cases are shown in Table~\ref{tab3}. In both cases, the prediction horizon contains a total of $48$ steps, therefore the problem size is identical. All Monte-Carlo parameters are as previously stated except for a reduction in the amount of mass variation due to the fragility of the fixed-timestep controller.

\begin{table}[t]
	\caption{Prediction Horizon Parameters} \label{tab3}
	\vspace{-1em}
	\begin{center}
		\begin{tabular}{|c|c|c|c|c|c|c|c|}
			\hline
			Case & $\Delta t_1$ (s) & $N_1$ & $\Delta t_2$ (s) & $N_2$ & $\Delta t_3$ (s) & $N_3$ & $T_\text{pred}$ (s) \\
			\hline
			1 & 0.04 & 24 & 0.16 & 12 & 0.64 & 12 & 10.56\\
			\hline
			2 & 0.02 & 48 & - & - & - & - & 0.96 \\
			\hline
		\end{tabular}
	\end{center}
	\vspace{-2.5em}
\end{table}

The distribution of RMSE results for $100$ Monte-Carlo runs is shown in Fig.~\ref{fig7}. With the longer total prediction horizon, the average attitude, velocity, and position errors in Case $1$ are $9\%$, $62\%$, and $64\%$ lower, respectively. However, the average thrust and torque input errors in Case $1$ are $81\%$ and $32\%$ higher, respectively. Comparing other metrics, the average time to reach the target and the average number of replanned trajectories is $44\%$ and $48\%$ lower respectively in Case $1$. Therefore, although the average control effort is higher in Case $1$, the longer prediction horizon achieved by the non-uniformly spaced timestep provides substantial benefit in tracking performance. Similar performance trends occur in no-disturbance simulations, but are not presented here due to space limitations.

\vspace{-0.5em}
\section{Conclusions}
\label{sec:conclusion}
An MPC approach for a tandem-rotor helicopter is presented in this paper. An augmented $SE_2(3)$ error definition is used to linearize the process model about a reference trajectory. A non-uniformly spaced prediction horizon is used to predict the multi-timescale dynamics while limiting the optimization problem size. The attitude is constrained using a combination of attitude keep-in zone and $\ell_1$-norm attitude error constraints. Monte-Carlo simulations demonstrate robustness to initial conditions, model uncertainty, and environmental disturbances. Additionally, the non-uniform prediction horizon is shown to be beneficial over a traditional fixed prediction horizon. Although a linear MPC formulation is used, the problem size still exceeds the limits of simple real-time computing platforms. Future work will focus on improving the disturbance estimation, enabling tracking of a moving target, and investigating explicit MPC approaches better suited for hardware implementation.

\bibliographystyle{IEEEtran}
\bibliography{IEEEabrv,lcss.bib}

\end{document}